\documentclass[11pt]{article}

\usepackage[letterpaper,margin=1.05in]{geometry}
\usepackage[T1]{fontenc}
\usepackage{libertinus}
\usepackage{microtype}
\usepackage{etoolbox}
\usepackage{setspace}
\usepackage{parskip}
\usepackage{enumitem}
\usepackage{xcolor}
\usepackage[most]{tcolorbox}
\usepackage{tikz}
\usetikzlibrary{arrows.meta,positioning,calc}
\usepackage{csquotes}
\usepackage[backend=biber,style=authoryear,sorting=nyt,maxcitenames=2]{biblatex}
\usepackage[colorlinks=true,linkcolor=linkblue,citecolor=linkblue,urlcolor=linkblue]{hyperref}

\definecolor{akcolor}{HTML}{244A64}
\definecolor{tkcolor}{HTML}{713B46}
\definecolor{linkblue}{HTML}{345D7E}
\definecolor{statementframe}{HTML}{687985}
\definecolor{statementbg}{HTML}{F4F6F7}
\definecolor{stagegray}{HTML}{5A6268}

\setlist{itemsep=0.35em,topsep=0.35em,leftmargin=2.2em}
\newcommand{\speaker}[1]{%
  \par\medskip\noindent%
  \ifstrequal{#1}{AK}{\textcolor{akcolor}{\textbf{\textsc{AK}:}}}{\textcolor{tkcolor}{\textbf{\textsc{TK}:}}}\enspace}

\newcommand{\stagebreak}[1]{%
  \par\vspace{1.2em}
  \begin{center}\textcolor{stagegray}{\small\bfseries\MakeUppercase{#1}}\end{center}
  \vspace{0.25em}}

\newtcolorbox{statementbox}[1]{
  enhanced,
  breakable,
  colback=statementbg,
  colframe=statementframe,
  boxrule=0.55pt,
  arc=1.5mm,
  left=2.2mm,right=2.2mm,top=1.3mm,bottom=1.5mm,
  before skip=0.8em,after skip=0.9em,
  fonttitle=\bfseries,
  title={#1}
}

\title{\vspace{-2.2em}\textbf{The Dialectical Theorem}\\[0.25em]
\large The future of universities, mathematics departments, and training of human mathematicians}
\author{\normalsize An informal conversation between\\[0.25em]
Thomas Koberda and Alex Kontorovich}
\date{}

\begin{document}
\maketitle
\vspace{-1.1em}

\speaker{AK} PhD students are hearing from all directions---from social media, from their peers, and from some famous professors---that ``math is solved'' and that mathematicians will no longer exist in the future (some even say humans won't exist in the future, but let's not derail our discussion by entertaining that notion; and anyway, if that is to be the case, then it's moot to argue about how many non-existing humans would be research mathematicians). So the PhD students hear: there's no point in finishing a degree; leave grad school now while you still can. And correspondingly, an undergrad considering grad school should think of something else to do. I want to think such arguments through in three stages:

\begin{enumerate}[label=\textbf{Stage \arabic*:},leftmargin=*,itemsep=0.45em]

\item (from the bottom up): Will math actually exist as a thing humans do for money? (My claim: Yes.)

\item (top down): Given that math expertise will exist, and hence that we will want to hire more people with that expertise, what is it that we'll be looking for when we hire?

\item Whatever the answer is to Stage 2, that expertise is what we should train PhD students to develop.

\end{enumerate}

[See  the end of the document for  a logical map of this conversation.]

\stagebreak{Stage 1}

\speaker{TK} I am profoundly worried about the future of mathematics as a profession, and I sometimes lie awake at night wondering whether, in a few years, there will be any jobs for human mathematicians, either in academia or outside it. I have friends, colleagues, and students who are currently on the job market; for them, these worries are particularly immediate.

Before we begin, I want to make clear—as I'm sure you agree—that this is not just an abstract philosophical exercise. I want us to articulate a point of view that I hope will become widespread, and that entails a future in which human mathematicians don't just exist, but are crucial. The thought experiment lies in the assumptions we allow ourselves to make, with the goal of arriving at conclusions that are as robust as possible.

With that in mind, I am happy, for the sake of this conversation, to take “math is solved” as an axiom.

\speaker{AK} Great, let's put down exactly what we mean:

\begin{statementbox}{The Solved Axiom}
Any problem that can be solved by human + AI can be solved by AI alone.
\end{statementbox}

I'm reminded of early experiments with ‘centaur’ chess tournaments, in which one human–computer team played against another. Pretty soon, computers became so strong that whenever a human departed from the computer's top suggested move, their team tended to lose. And obviously that's no fun, so we don't see many such tournaments anymore.

Now, I don't actually believe that the Solved Axiom is true today, nor do I believe it will be in the foreseeable future! But I feel compelled to grant it for the sake of argument. Even {\it with} this axiom, I  still want to argue that math is not dead and that human expert mathematicians will continue to be in demand in society.

In fact, let's give the premise I'm arguing against a name:

\begin{statementbox}{The Doom Premise}
Assume the Solved Axiom. Then there is no longer any demand for professional human research mathematicians or human mathematical expertise.
\end{statementbox}

There is a related prediction that people often seem to make, and I think it is worth stating separately.

\begin{statementbox}{The Training Premise}
Assume the Solved Axiom. There is no longer any demand for training in mathematics for PhDs, undergraduate majors, and so forth.
\end{statementbox}

I should of course admit that I neither believe the Doom Premise, nor do I {\it want} it to be true, and so, to make my argument against it as strong as possible (and to protect myself, as much as possible, from wishful thinking and other biases), I'll need to argue purely logically and dispassionately, granting as much as I can to the ``Devil's Advocate'' and relying only on pragmatism: not what I {\it wish} would happen, but what I see as most likely to happen when the dynamics of voluntary action play out.
So, I think the strongest version of my argument against the Doom Premise grants the Solved Axiom (which is anyway an assumption in the premise), and continues from there.

\speaker{TK} Even if only a weaker version of the Doom Premise ends up being true—for instance, a dramatic reduction in federal funding for mathematics—the consequences could be severe. But I agree that we should test our argument against the strongest version of the Doom Premise, for the sake of robustness. I would also like to isolate one further possibility, concerning the economic support of mathematics, and make it explicit.

\begin{statementbox}{The Funding Premise}
Assume the Solved Axiom. Public agencies, philanthropic organizations, and corporations will no longer fund mathematical research conducted by humans.\end{statementbox}

We will need to talk about the Funding Premise more, but let's address that a bit later.
For the sake of the conversation, we can expand the Solved Axiom to mean that the AI is much better at explaining to humans what it has done than any human interlocutor, that it is better at posing the next natural directions of research, and that it guides its own research program much better and faster than any team of humans could.

\speaker{AK} I would make a distinction between the production of new raw knowledge and the capacity to teach humans. (In fact, I will later argue that humans, in the main, do not want to sit for hours a day staring at a screen trying to learn, as enticing as those screens may be... [TK: Let's see if we agree or disagree about this when we get to this topic!]) But sure, yes, let's say that AI is also a better teacher (of anyone, including humans), as well as more knowledgeable and more patient than any human. (This is probably already true today! It's just that I'll argue later that people will not want to work with AI the same way they'd work with a human tutor).

\speaker{TK} I agree that the Solved Axiom does not describe the present state of AI, and it is unclear whether current AI architectures could ever satisfy it, or whether something fundamentally different would be required. Nevertheless, I don't want our argument to depend on predictions about the future development of AI. I would much rather arrive at conclusions that remain valid even if the strongest version of the Solved Axiom eventually becomes reality.

\speaker{AK} Great, yes, let's accept it in all of its forms. Then my first claim---and I will want to hear your pushback once I've had a chance to argue for it---is:

\begin{statementbox}{Expertise Demand Theorem}
Assume the Solved Axiom. Even in that world, there will still be demand for human mathematical expertise.
\end{statementbox}

This, of course, is a direct negation of the Doom Premise.
If we can successfully establish this theorem and hence human experts are to continue to exist, well, they'll need to come from {\it somewhere}. So we'll follow that up with the corollary:

\begin{statementbox}{Training Demand Theorem}
Assume the Solved Axiom. Humans will want to be trained to become experts in mathematics.
\end{statementbox}

And this, of course, is the direct negation of the Training Premise.

\speaker{TK} I'll be very interested to see how these theorems will be established. I think my greatest worry for now is that there will be a collapse, or at least a rupture, in the process of training mathematicians, after which recovering a culture of humans doing mathematics will be impossible. We could see something like this play out: universities broadly decide that mathematics departments and mathematicians are unnecessary, both as researchers and as teachers; funding agencies decide that mathematics research is not worthy of societal support; and loud and powerful voices in tech decide that mathematicians (and indeed most human intellectual activities) are obsolete.

Have you seen that the University of Chicago considered (though, it seems, did not implement) plans to eliminate large parts of foreign language instruction and replace it with ChatGPT? See \textcite{ando2025crisis}.

\speaker{AK} Some universities may well make such choices! (But why in the world would I pay tuition to U of C just to sit on GPT and learn French?) I'm not against anyone making such decisions (of course I actually am! But I don't need to be for the sake of the argument I'm trying to make); everyone should make their own decisions and face the consequences. My feeling, which I will try to argue and elaborate on below, is that universities that choose not to downsize their mathematics offerings will, over time, do better than those that do. As long as not everyone at the same time makes the (terrible, in my opinion) decision to close down their math departments, and as long as some professional expertise still exists, then the value of that expertise will have a chance to prove itself in the marketplace of ideas and practice.

\speaker{TK} I think a lot of work is being done by ``not everyone at the same time.'' I think it is incumbent on us as a community to have a compelling argument to counter the inference ``math is solved, therefore human mathematicians are unnecessary,'' both to prevent everyone from making the same terrible decision and to ensure that the notion that mathematical training is not just useful but a good for both individuals and society remains widespread.

\speaker{AK} I agree completely that if all mathematicians simultaneously ceased to exist, it would take a long time and a lot of work before (if ever) we could regain mathematical culture and practice. I just don't believe this will happen unless many people all make the same bad decisions at the same time. An analogy is to arguments 80 years ago about whether children need to learn their times tables, given that mechanical calculators can now do arithmetic. Some educators did---and some still do---argue that children need not memorize them. Those children and their teachers listened, made bad decisions, and suffered the consequences. For any children under my care, I insist that they know their times tables, and I (hopefully) then watch them thrive mathematically. The outcomes of such decisions are not seen immediately; they play out over years and decades. So here, too, it will take time for people to see the consequences. If some universities choose to close down their math departments, they might not immediately fail; but I believe (and will argue more precisely when we get to it) that those universities that choose to invest in their math departments will indeed do better, and hence professional math expertise will continue to exist.

\speaker{TK} A decades-long process seems to me like it could be enough time to cause a full-scale break in the continuity of mathematical culture, from which recovery would be extremely difficult, if not impossible.

\speaker{AK} I agree that if many people all independently make very bad (from my viewpoint) decisions, then this is a possibility. I just don't think there's a high likelihood of so many bad decisions being made. Some people will ``hold out,'' and we'll be able to see that those who do excel, and so others may rethink their decision. There's nothing I can do to counter such decision-making except to speak publicly about why I think 
math should continue to be supported.

\speaker{TK} Great, so let's figure out how to convince people. One immediate complication I see is that adults deciding whether or not to study higher-level mathematics are not children, and they will inevitably be influenced by many more voices and opinions, some giving good advice and some giving bad advice.

\speaker{AK} I agree! The important thing is for there not to be just a few people charged with making such decisions for all of us. Each university president would need to make the decision independently, not follow some edict from the Department of Education (or whatever). Even if all NSF funding disappears, that doesn't fundamentally change my argument. In fact, I'm willing to grant the full Funding Premise: suppose all external funding for human mathematical research disappears, including support from industry and philanthropy. Sure, let's assume it! I still think we'll be fine. [TK: I'm less certain that we will be fine, but I'm willing to consider that possibility for the sake of the conversation.] Of course, I would much prefer that such funding continue, both for my own selfish sake and because I believe it is actually valuable for society. \cite{etingof2026} describe human mathematicians as an ‘intellectual reserve,’ which I think is an excellent way of expressing at least one reason why preserving human mathematical expertise matters.

\speaker{TK} I also like the idea of intellectual reserve. Building on that, one point that I made in my essay about graduate training \parencite{koberda2026training} is that retaining human experts in mathematics is crucial for preserving human epistemic self-governance, which I suppose could be summarized as follows: we need to understand math if we want to have any idea of what is going on in the world around us. I am curious what you think about this and whether it could provide a basis for an argument that would convince society that training human mathematicians is a valuable activity.

\speaker{AK} I guess where we differ is that I'm not---in the argument I have in mind---trying to convince ``society'' of the importance of training humans. Not because I don't think it's important (I think it really is! And I do everything I can to do it!), but because I don't think it's necessary for my argument. I think that, as long as individuals are capable of making choices and seeing the outcomes of those choices, there will be people who choose to study mathematics, which means they will want to learn from expert human mathematicians, which means there will be demand for such experts.

\speaker{TK} Let's expand the point about learning from expert human mathematicians first. You said earlier in this conversation that you believe that humans will definitely prefer learning from humans rather than from AI/computers. Let's clarify this.

\speaker{AK} Yes, let's make this an explicit axiom:

\begin{statementbox}{Learning Axiom}
On the whole, human beings who want to learn deep mathematics will want to learn it from other humans. 
\end{statementbox}

\speaker{TK} We are mostly talking about mathematics in this conversation, and so we should keep our conversation focused, though I would not argue that this axiom is specific to mathematics.

\speaker{AK} Yes. But first: the Learning Axiom does not mean that there won't be ``Ramanujans'' in India who iterate with an AI at their public library and deeply learn automorphic representations; in fact, I do believe that such things will now be possible, and that AI has the potential (though this is certainly not guaranteed!) to democratize access to the highest levels of knowledge for people who are geographically, socially, or otherwise far removed from experts. But my claim is that {\it most} people will not be able to do this; we are basically furless apes, and we want to be around other humans, talk to humans in person, learn from human experts, watch our human peers struggle and overcome, and enjoy all the social things that come with that. My biggest evidence for this axiom is COVID: we put not AI but actual live humans on Zoom screens in front of millions of children, and the experience was, on the whole, enormously unpopular! See also [as suggested by Chris Connell] the flop of MOOCs [Massive Open Online Courses].

\speaker{TK} My worry about accepting the Learning Axiom is not that it contradicts the observation that students (from elementary school on up) hated learning even from other humans on a computer screen, much less from a chatbot. Rather, under an axiom as broad as the Solved Axiom, why would one rule out the possibility of an AGI that has ``hacked'' the social circuitry of the human brain and is a much better conversationalist than any human could hope to be, to the point that even the most AI-averse human would concede that they feel heard and understood by the machine better than by most, if not all, of the humans they know?

\speaker{AK} I'm perfectly happy to grant that: for humans who have the self-drive and desire to iterate with an AI, they will learn better than humans who learn from social contexts. I just don't believe that most humans have that self-drive, and most would be much happier (and thus would voluntarily choose) to be around other humans. That's what the axiom is trying to capture.

\speaker{TK} I do think I have to disagree here, though. I am not convinced that a pedagogical culture prioritizing interaction with other humans will survive automatically; I think it may need to be protected with guardrails. We can probably agree that completely unsocialized humans who only rarely interact with other humans are undesirable, but I do not think we can simply assume that such a future will not happen.

\speaker{AK} Again, what I'm really trying to do is not argue for the things that I {\it hope} will happen. I do hope that my conclusions will happen, but my motivation is instead just to follow, to the best of my limited ability, logical deductions from the most basic axioms I can think of, wherever they may lead. Some of the axioms I don't actually believe but am willing to grant (e.g., the Solved Axiom); others I very much believe to be true but don't know a more basic reason for them to be true, and they're necessary for my argument, so I am forced to take them as truths I hold to be self-evident (e.g., the Learning Axiom).

As to whether human interactions should be prioritized with guardrails: I disagree. I think people will just want these things, because it's our nature. We shouldn't need any formal guardrails; nature has given them to us.

\speaker{TK} I hope you're right and that people will realize that socialization through interaction with humans is necessary for human flourishing, but I still retain my caution here.

\speaker{AK} As with other dynamical systems, I think some people will choose to sit in their rooms alone, iterating with AI. They're welcome to do that, even if it may be to their own social detriment!  I just don't think they'll be as happy as people who are socializing. And so I am claiming that people will voluntarily (again, for the most part) choose to participate in experiences where socializing is plentiful.

\speaker{TK} Great, so with the Solved Axiom and the Learning Axiom, why do you think that universities will continue to exist?

\speaker{AK} Okay, I think I need one more axiom.

\begin{statementbox}{Value Axiom}
People who are capable of deep, extended thought will, on the whole, have richer lives than people who are not.
\end{statementbox}

This does not necessarily mean that they're economically richer (I do believe that as well, but don't need it for my argument). [TK: So maybe more intellectually and emotionally fulfilled? AK: Yes; in terms of general human flourishing.] What I actually believe is much stronger: people who are capable of deep thought will be much better off (in the main) due to the existence of technology that will make their skill even more valuable. But again, I don't need this stronger form for my argument.

\speaker{TK} And so, from the Value Axiom and the Learning Axiom, you want to infer that universities that train human intellectuals (and mathematicians in particular) will continue to exist.

\speaker{AK} Yes, exactly! I suppose I may need another axiom: that one of the best ways to train one's mind to be capable of deep, sustained thought at this highest level is the training one experiences (or should!) in a math PhD program.

Do you not want to argue against the Value Axiom? In some sense, it goes counter to the Solved Axiom. Doesn't the latter mean that humans have no value to add? Hopefully the general sense of the words ``richer lives'' is doing the bulk of the work. Does that make sense?

\speaker{TK} Yes, I do think that there is a tension between the Value Axiom and the Solved Axiom, though how much depends on how broadly we interpret the semantics of these axioms. In a solved world, one could suggest that there are other ways of cultivating deep, extended thought and achieving the attendant intellectual and emotional fulfillment. Why would this have to be through mathematics at all, especially if all the intellectual paths are well-trodden and completely understood by the AI? A lot of the friction necessary for developing mathematical skill is unpleasant, inefficient, and sometimes even a waste of time. Why would we not look somewhere else to enrich their lives? Without implicitly acknowledging the necessity of something like epistemic self-governance, how do we argue that the Value Axiom and the Solved Axiom are compatible?

\speaker{AK} Okay, so let's make the other axiom explicit.

\begin{statementbox}{Training Axiom}
One of the best ways to cultivate the ability to think deeply and clearly is through the study of pure mathematics, the more abstract the better.
\end{statementbox}

The evidence for this includes the fact that Abe Lincoln apparently studied Euclid to sharpen his rhetorical skills\footnote{From \parencite{nytimes1864}: ``At last I said [to myself], `Lincoln, you can never make a lawyer if you do not know what {\it demonstrate} means;' and I left my station in Springfield, went home to my father's house, and staid there till I could give any propositions in the six books of Euclid at sight. I then found out what `demonstrate' means, and went back to my law studies.''}. I agree completely that math isn't the only path to clear thinking, but it's a pretty darn good one. Even if writing down solutions to mathematical problems need not be difficult anymore (because AI can do it), students who choose to do the exercises as intended, with all the friction, frustration, confusion, etc., are those who will actually achieve the effects of the training. (This will all be discussed in detail in Stage 3.)

\speaker{TK} So, as I understand it, your argument is fundamentally pragmatic. [AK: Yes!] There are perhaps many ways to achieve the ability to think deeply (in my essay I put this in a somewhat broader category that I called ``intellectual agency,'' in order to distinguish it from ``independence,'' for instance), but mathematics already exists as a path to this goal, and so we might as well choose that path.

\speaker{AK} I agree completely. First we have to establish that the ability to think clearly is a goal (Value Axiom), then that learning mathematics is an efficient mechanism for attaining that goal (Training Axiom), and finally that people will want to undertake this training in person (Learning Axiom).

\speaker{TK} I am still uncomfortable without a stronger justification, but, more seriously, I worry that this does not fully justify why the Value Axiom and the Solved Axiom are compatible. More specifically, how do we justify the goal? I feel like there is some implicit acceptance of the necessity of epistemic self-governance on the part of humans.

\speaker{AK} You're right on both counts. I don't know how to justify the Value Axiom from more basic principles, hence its status as an Axiom. Do you think you're willing to grant it? (In terms of human self-governance, I do hope that we continue to live in an open society where people are, largely, free to choose what they want to do with their time. Such questions are really important, but I think they are outside the scope of our discussion.)

\speaker{TK} I would like to grant it, but I think that the following axiom would be more fundamental, and I would welcome your disagreement:

\begin{statementbox}{Epistemic Axiom}
Humans must retain deep understanding in important domains of knowledge and remain capable of functioning as independent intellectual agents; this is not to say that humans must be able to function without AI, or that they should not use it. More simply, humans should be able to think for themselves sufficiently well that they can understand and evaluate non-human intellectual activity.
\end{statementbox}

I would say that the Value Axiom could be a corollary of this one. Moreover, it would seem to justify accepting the Value Axiom under the regime of the Solved Axiom.

\speaker{AK} If you see the Value Axiom as a corollary of the Epistemic Axiom, then I'm perfectly happy to grant the latter and get access to the former! I don't believe that there is any value in a massive database of mathematical theorems and proofs---even correct, well-written, and formally verified ones---that exists only as 0s and 1s inside data centers, but not in any human heads. That is, using AI to create new mathematics for its own sake, without having any human understand it, will not be something people do, certainly not once the novelty of it wears off. A  point of comparison: millions of kids around the world will dress up on Halloween as Harry Potter because they've read the books, watched the movies, and been deeply moved by the stories, which have brought great joy to children (and adults!) around the world. It is easy to imagine, especially in light of the Solved Axiom, that AI could write thousands of stories that are every bit as moving. If we let it write those stories but no human ever reads and enjoys them, what was the point? The required energy expenditure would be a complete waste; who would pay the electric bill? Similarly, I claim that mathematics is only valuable insofar as it exists in the brains of humans. Reading Flexner's ``The Usefulness of Useless Knowledge'' \parencite{flexner1939} (written in 1939 to explain why someone might want to support a ``useless''-knowledge-generating institution like IAS), I think that all of his arguments fall apart without changing the title to ``The Usefulness of Useless Knowledge {\it in the Minds of Humans}.'' Who could have known in 1939 that this add-on would become necessary!

\speaker{TK} I think that we mostly agree here, though I might soften the formulation slightly. It was already true before AI that mathematics had to be alive in order to realize its value, which is to say that it had to form part of human mathematical life. A theorem sitting unread in a library may still have potential value---it could be rediscovered decades later---but that value is realized only when the mathematics is understood, used, taught, extended, or otherwise brought back into human intellectual activity. Here I would emphasize that the Epistemic Axiom is playing a key role, since otherwise one could argue that mathematics sitting in some data center is still useful to an AI that then uses it to develop some other new technology that benefits humanity (or possibly serves the AI's own ends).

Historically, areas of mathematics have gone in and out of fashion, and the skills and knowledge possessed by practicing mathematicians have evolved over time. Inevitably, some knowledge has been lost. This is one reason it is so important to write clear papers, and why exposition that organizes bodies of literature and makes them usable to researchers approaching from the outside is so valuable. It is a shame when mathematical writing languishes indefinitely in libraries or preprint repositories. Part of the work of mathematics is to keep such knowledge alive by absorbing it, rediscovering it, and teaching it to the next generation of mathematicians.

I think this gives us a great segue: how do we now infer the future necessity of universities? I do agree that universities provide a lot of the infrastructure (at least ideally) that is conducive to training, including fostering intellectual culture, enjoying relative independence from commercial pressures, creating opportunities for socialization, and providing long time horizons. It might be the case that, in the future, competing institutions will arise that do these things as well or better, and that would be welcome, but this might be more of an argument for not immediately abolishing universities than for the more positive conclusion that ``universities will persist in the future.''

\speaker{AK} For my purposes, the very definition of a ``university'' is a place that ``provides a lot of the infrastructure that is conducive to training, including intellectual culture, relative independence from commercial pressures, socialization, and long time horizons,'' and a place where experts congregate for the purpose of producing new (and retaining old) deep knowledge and disseminating that knowledge to students eager to learn it. To the extent that existing institutions currently called universities do these jobs, I think there will continue to be demand for them. If they don't do these jobs, or if new institutions arise that do them better and students choose to go there instead (I'm all for it, for the sake of argument), then {\it which} institutions exist will change, but what they {\it do} will not. [TK: I agree with this.] So the main point of this line of reasoning is to finally prove the Expertise Demand Theorem, and hence its corollary, the Training Demand Theorem (which, now that I think of it, provides a positive feedback loop into the Value Axiom + Training Axiom). In fact, my personal opinion is not just that universities with mathematics departments will continue to exist, but that there will be even {\it more} demand for mathematicians in the future than there is today, precisely because the training by such experts will be so valuable. But I don't want to add this stronger argument to make my main points here.

Are we ready for Stage 2? Given that universities (maybe not all the ones that currently exist, and maybe some new ones that haven't been invented yet) will need to continually replenish their collection of experts, what exactly would hiring committees be looking for? What do we, in the world of mathematical culture and practice, value? (That, of course, will feed into Stage 3: what one should be trained in during a PhD.)

\stagebreak{Stage 2}

\speaker{TK} Yes, I think I am satisfied with our answer to why universities (as defined above) will continue to exist. Now, as to why they would continue to hire mathematicians and what will be valued in hiring: at least the first of these seems to have an answer that follows from what we have said. If we accept that there will be demand for humans to be trained in mathematics by human experts, then universities will presumably hire human experts in mathematics simply to meet that demand. What is less clear to me is which criteria will be used to evaluate and hire mathematicians, which is precisely the values question. I think this ties in very closely with the question of what the purpose of graduate education should be, since presumably graduate education should function to train human experts. The nuances of the function of graduate education are numerous, and I made some effort to clarify them in my essay. It seems safe to assume the following:

\begin{enumerate}[label=\arabic*.,leftmargin=*]

\item A mathematician hired, retained, and promoted by a university will be expected to have intellectual agency, as defined and elaborated in my essay.

\item The mathematician will be expected to have contributed to human understanding of mathematics. This is a broad category, which would include traditional metrics such as original papers. Part of the issue at hand, of course, is that the existence of papers authored by the mathematician does not serve as the reliable proxy for intellectual agency that it did until recently, which is one of the reasons I put intellectual agency first.

\item A mathematician will be expected to have excellent skills in communicating mathematics to others and continued evidence of contributing to human understanding of mathematics. Communication of mathematics can include mathematics created by AI, but the human mathematician should play a critical role in translation.

\end{enumerate}

\speaker{AK} I would prefer not to assume these things, that is, not make them Axioms, but rather conclude them as Theorems from more basic premises.

\begin{itemize}[leftmargin=*]

\item For point 1, that's part of our working definition of what a university is, so I think it does not need to be a separate assumption. (In Lean-speak, the proof is \texttt{rfl}, that is, by reflexivity.)

\item For point 2, we can argue by the contrapositive: what would it mean if mathematicians were not expected to contribute to human understanding of mathematics? In light of the Solved Axiom, this presumably means that the so-called ``expert'' (perhaps we should just switch to using the standard word for it, Professor?) has pushed some buttons in GPT42, produced a slew of papers solving longstanding problems, and put them on their CV, while having no clear, deep understanding of the mathematics those papers claim to establish. Why would a hiring committee make an offer to such a person? I can press a button in GPT, too! (Remember when ``Prompt Engineer'' was a career, for a hot minute? Now amateurs type ``do a brkthru'' [internet subculture slang for prompting AI to prove a breakthrough result] and actually get one! No thanks.) It usually takes around 3-5 years of hard work to first internalize some deep piece of difficult mathematics. It takes 3-5 minutes of a job talk to determine that someone has not done this work.

\item For point 3, if the reason the university exists is that students are paying tuition to go and learn (cf. Learning Axiom), then a Professor who does not communicate mathematics well is not useful or valuable to the enterprise. (Moreover, how could someone communicate effectively with critically thinking students who can ask penetrating questions at any time if they do not themselves understand the material deeply? This reinforces point 2.)

\end{itemize}

\speaker{TK} Agreed! This is why I am saying it should be safe to make these assumptions; maybe infer is better than assume. [AK: Okay, great, then we're fully in agreement.] Why don't we call these three points the Hiring Value Theorem? Of course, as hiring seasons proceed in succession, we will eventually gain an empirical demonstration or refutation of the Hiring Value Theorem, since if hiring does actually continue in the future, some people will secure positions, and the community will be able to infer {\it ex post facto} what was evaluated.

\speaker{AK} Yes, let's state the theorem explicitly.

\begin{statementbox}{Hiring Value Theorem}
In order to be hired, retained, and promoted by a department of mathematics, a person will need (1) to understand deep, difficult mathematics, and a lot of it, and (2) be able to communicate that understanding effectively to others, including verbally!
\end{statementbox}

The oral part is the new focus of attention. For thousands of years, if you could communicate in written form, that was a nearly perfect proxy for deep understanding. (I suppose you could've paid someone to write something for you, then crossed out their name and put your own. I'm looking at you, l'Hôpital! It's Bernoulli's Rule, not yours!) Moreover, almost the opposite was true: if you were a brilliant communicator of your ideas but failed to document them in written form, you would not be contributing to the scientific literature and hence could not get a job! Since AI systems can (in the world of the Solved Axiom) write amazing, perfectly well-written texts, the correlation between producing written arguments and understanding has plummeted to zero. The written argument may remain every bit as mathematically valuable, but what has changed is its value as evidence about the person whose name appears on it. One particularly effective way to demonstrate, in practice, that you really understand a piece of mathematics deeply is to give a lecture on it to experts (and, of course, also publish a corresponding paper in a good journal; I don't think we want to discard the scientific literature!) and adequately field all of their most penetrating questions. This is a really important point, I think, so let's give it a name:

\begin{statementbox}{Dialectical Understanding Axiom}
Assume the Solved Axiom. Then one of the principal ways a person can effectively demonstrate to another that they understand a piece of mathematics deeply is to explain it orally, and respond coherently to questioning.
\end{statementbox}

This certainly does not mean hostile questioning. One of the things I love most about math culture is how open we are to this kind of dialectical questioning (just like what you and I are doing right now!). In fact, I'm very careful when early graduate students are exposed to it for the first time: it takes a lot of experience to get  used to hearing someone say, ``No, what you've said is wrong, and let me try to argue why,'' and being not only perfectly fine with it but actually valuing it (because it means you're about to learn something new). Indeed, it's radically different from the norms of conversation outside mathematics---as I'm painfully reminded from time to time whenever I instinctively respond to a non-mathematician in an unintentionally rude-sounding way. Our most precious expression is even called {\it Proof} (or {\it Argument}), so we're supposed to have that dialectical process happen internally before it even gets to the world outside ourselves!

\speaker{TK} I feel that at this point it is important to note that there are dangers lurking here, especially in how the Hiring Value Theorem is applied in practice. This has to do with the inevitable shortage of positions relative to the number of qualified candidates, and so there will inevitably be questions about what is judged more deserving of attention from human mathematicians. Relatedly, there is a hazard of exacerbating inequalities that the community in general has already fought hard to correct. There is a danger that a small handful of people could become authoritative arbiters of a nebulous and opaque measure of importance that we might call ``taste.'' If we evaluate job candidates based on oral performance only [AK: Or even mainly.], we need to make sure we do not simply reward the eloquence of native-language speakers, for instance. These are all important and pressing questions that I want to flag, though I don't think they are the central concern of this discussion.

\speaker{AK} I agree completely. Oral evaluation inevitably introduces various irrelevant human advantages—eloquence, confidence, linguistic background, attractiveness, social ease—and I don't like that. But I don't see how we could eliminate them completely. We should nevertheless work as hard as we can to evaluate deep mathematical knowledge and communication on their merits.

\speaker{TK} Yes. Irrelevant human advantages make me uncomfortable, but I also do not have a good fix beyond vague suggestions like ``transparent standards'' and ``multiple evaluators.'' Working out those details is an important task for the mathematical community.

\speaker{AK} Honestly (and this is perhaps getting us off track), I also don't see a way to ``fix'' this. Those who choose to abuse the system and select more for, say, ``height'' or the ``right regional accent'' (especially in the U.S.) instead of ``mathematical knowledge and communication'' will, over time, have a less highly ranked department. So the dynamics in the system should already incentivize meritorious hiring.

\speaker{TK} Should we leave this part of the conversation for another day and move on to Stage 3?

\stagebreak{Stage 3}

\speaker{AK} Yes, great, let's! So finally we've recursed far enough back to reach the point of this conversation! What is the purpose of a PhD education? Certainly one of its purposes should be to prepare someone for the possibility of an academic career (simply because there's almost no other way into academia). (Of course, even pre-AI, the vast, vast majority of people who received a math PhD did not go on to a tenured professorship. Those who left academia often found their skills in deep, abstract thought very valuable, gaining meaningful employment, making much more money than professors, and being able to choose where they want to live; we silly professors have the funny habit of calling this ``failing'' to get an academic job! [TK: It's part of our cult-like academic culture!] Anyway, now I'm rehashing the arguments for the Value Axiom, which we've already agreed to assume.)

Given that this is the case, how does one train a student in deep mathematical knowledge and oral communication (since that's what hiring committees will be looking for)? I think the answer is pretty simple: not all that much needs to change! What we value is deep understanding of difficult mathematics, and acquiring that understanding still takes years of really hard work, with AI help or without it. What may change radically is how we determine whether someone has actually acquired that understanding.

Pre-AI, a typical five-year PhD student experience was something like this: two years of preparatory work, coursework, reading books, and doing written and oral exams; followed by two years of unadulterated research on a problem (or a few such), starting to go to conferences and learning the field---what questions people are interested in, and so on; and then the last nine months were spent writing up a proof in complete detail, getting everything down on paper, learning how to organize an argument, and putting the finishing touches on a thesis; and finally, defending the thesis, which is basically a nice, easy talk on your work (which your advisor, if they know what they're doing, has already certified as sufficient to warrant a PhD).

So the written thesis was a nearly perfect proxy for this: the student had done a significant amount of hard work, learned deep mathematics, contributed something new to our collective understanding, and struggled to prepare a carefully written document. Under the Solved Axiom, you can now get such a document with the push of a button. So how should one prove that they've learned a great deal of mathematics? We've already covered this: the Dialectical Understanding Axiom!

This means we can finally state the main argument I was hoping to make:

\begin{statementbox}{The Dialectical Theorem}
What you should do in grad school is learn a ton of math and, as much as you can, expose yourself to good-faith questioning to ensure that you know your subject inside and out and can effectively communicate it to others.
\end{statementbox}

Later I'll give some more concrete suggestions for how to do this.

\speaker{TK} I think you make good points here. Thinking about difficult problems that resist many attempts at solution creates a type of friction that requires persistence, patience, and personal growth to overcome. Of course, personal growth does not mean that everything has to happen in isolation; indeed, it might be even better if the process is dialectical. The thesis was a reliable certificate that this sort of friction was encountered and the corresponding growth occurred. It seems that some sort of cognitive friction and engagement with difficult problems is the principal path to achieving what we are calling intellectual agency, and it is unclear to me whether there are effective alternative paths, or whether any amount of simplification, easing, or expediting of the process could be achieved with AI. [AK: Agreed! Maybe in the future, we'll have Matrix-style brain plugs; but until then, I also don't see alternate paths!] One idea gaining currency is that there are two main modes of training that one could implement for students, namely independent [AI-free] and augmented [with AI allowed] modes. Such a formalized distinction would give a framework for distinguishing skills that need to be acquired through cognitive friction from those that can safely be offloaded at a later time, so that AI does not simply become a veil obscuring a lack of deep understanding.

\speaker{AK} I like that! As for how to use technology (AI) well, I like to tell my students the following three vignettes:

\begin{enumerate}[label=\arabic*),leftmargin=*]

\item Say you want to move some heavy pallets from a warehouse onto a truck. You're not going to do that by hand. So you grab a forklift and drive it over, pick up the pallets, and load up the truck.

\item Say you want to learn how to use a forklift. So you get in your forklift and drive over to some pallets in a warehouse, pick them up, and move them to another part of the warehouse, not because you wanted to move them, but because it was good practice in the use of that technology.

\item Say you want to get a good workout and get stronger. So you go to the gym and load 200 lbs. on the bench press. Then, instead of lifting the bar yourself, you grab your forklift, drive it over, and use it to do 20 reps of lifting the bar.

\end{enumerate}

I hope you agree that 1) and 2) are perfectly fine uses of technology, whereas 3) is extremely silly and counterproductive. To me, solving a hard open math problem, creating new mathematics for your thesis, etc., is akin to vignette 1. Learning how to prompt AI systems is Type 2 (which, by the way, I'm rather opposed to PhD programs deciding to teach---students will pick up prompting on their own; it's really not that hard!---and the precious time spent in grad courses should be reserved for learning deep mathematics). Doing homework problems in a course by prompting an LLM for the solutions is Type 3. If you're using AI to answer homework problems, then ask yourself what the point is. That said, even in Type 3, you could choose to use technology meaningfully: in the gym, you can ask AI to help you work out a training regimen, ask it questions about good diets, etc., etc. In the same way, you can ask an AI to clarify something the professor said in class that you didn't understand, give you more practice problems, etc., etc. Just don't use it to bypass the difficult, friction-filled stuff!

\speaker{TK} Yes, the vignettes make sense, though I'd like to elaborate a bit more, especially on Type 3. I think a stronger analogy with the goals of mathematical training would be that the task at hand is to get physically stronger, stay strong, and move heavy pallets laden with gym weights from a warehouse, all at the same time. To get stronger, you need to do physical exercise so that your body adapts. Once you have gained strength, you will want to maintain that physical strength, but that does not mean that you have to transfer a warehouse full of pallets by hand. Moreover, individual variations in desired levels of strength, needs for maintaining health, etc., would inform when it is better to do lifting by hand and when it is better to do it by forklift. This highlights the individuality of what should be considered ``premature cognitive offloading'' versus ``productive cognitive offloading.''

I would probably add that graduate students need to get to know themselves and collaborate with their PhD advisors to learn what constitutes premature and productive cognitive offloading when it comes to AI usage.

\speaker{AK} In my experience, it's very easy in practice to tell whether what you're doing is premature or productive. There's an angel on one shoulder and a devil on the other, and they both speak very clearly: one says, ``Hey, you shouldn't be asking for this answer; you know you should try to figure it out yourself!'' and the other says, ``That was a good practice problem for us to iron out; good job getting the AI to ask it!''

\speaker{TK} I'll take your word for it. [AK: Ha!] As for recognizing this in yourself, I might suggest that this is a learned skill that follows from already being highly trained as a mathematician and that it need not be automatic.

\speaker{AK} Fair enough, I'm happy to concede that! Actually, with my grad students, I can quickly tell when they've gone off, pushed a button in AI, written down what it said, thought about it a little, but didn't really try to argue all sides of the claim. So when they come to see me the following week and I start asking basic questions, they often turn out to be questions the student hasn't considered, but should have if they had interrogated their own arguments fully. To be fair, this was also sometimes an issue pre-AI, but it's perhaps even more prevalent now.

\speaker{TK} For now, at least, it seems that this is one of the concrete benefits of having humans train other humans: they can teach graduate students to interact more critically with AI. [AK: Agreed!]

So what do you think university mathematics departments will look like in the future? Will there be fewer jobs with even fiercer competition? Or do you think universities will be incentivized to hire even more mathematicians than they do currently?

\speaker{AK} Well, I've tried to keep this discussion to basic things I can derive from what I hope are well-formed Axioms, and to conclude only the weakest form of ``the world [mathematical or otherwise] will not end.'' But since you're asking, I'm personally extremely optimistic about what math will look like in the future. And as I see it, unless we all make a long series of very bad decisions, independently and one after the other (which I find highly unlikely), we are headed for a really enjoyable world of plentiful mathematics, with lots of work to be done to digest, understand, and synthesize beautiful arguments. [TK: Which could be either AI- or human-created, of course.] So I would hope that we'll soon see much greater demand for many more mathematicians than we currently have! What do you think?

\speaker{TK} I think that what you are saying is that if AI dramatically increases the amount of mathematics humans can produce, understand, and use, then demand for human mathematical expertise need not shrink with the amount of labor required for any one task; it might even grow. [AK: I suppose that depends on the demand function, but yes, I do think something like that will hold.] I hope that this is valid, because this is the most optimistic view of the future of which I can conceive. I think our conversation has made a compelling case that human mathematical understanding will remain valuable, even under the strongest assumptions about AI. What I remain less certain (and even quite pessimistic) about is whether the institutions and economic incentives necessary to sustain a thriving mathematical profession will endure. My worry, as I've articulated earlier, is that too many universities and individuals will make decisions that, in the aggregate, become very difficult to reverse. This discussion is a small part of a much larger story that is playing out and whose conclusion is far from certain; see for instance~\cite{AbbesMiller2026}. That is partly why I think we need to have this conversation now, rather than waiting for the consequences of those decisions to become clear. I hope our conclusions stand the test of time.

\speaker{AK} Me too! There's one more thing we haven't quite done yet. We've talked at length about what graduate training should accomplish, but perhaps we should finish with a few concrete suggestions for students who are wondering what they should actually do now.

One thing they can do right now is to make a habit of explaining the mathematics they're learning to other people. Take a theorem you're working through, explain it at the blackboard without notes, and encourage your fellow students to ask difficult questions. Why is that hypothesis necessary? What happens in particular examples? Where might the proof break down? These are things you can practice long before you have a thesis to defend.

\speaker{TK} I would add that students should deliberately make time to work independently, even when AI assistance is available. One good exercise might be to return to an argument that an AI helped you understand a few days earlier and try to reconstruct it without assistance. If you cannot, that tells you something important about what you still need to learn. Of course, how much independent work is appropriate will vary between students and between different stages of their training.

\speaker{AK} And I would encourage students to participate as much as possible in the mathematical life around them: attend seminars, attend the weekly department colloquium, attend conferences and workshops, ask questions, argue about proofs with their peers, and learn mathematics well enough to explain it to someone else. These things have always been part of becoming a mathematician, and I think they will be more important than ever.

\speaker{TK} I agree. And I think we should emphasize that advisors and departments have a responsibility to create environments in which students can do these things. Especially early in their training, students need opportunities to ask questions, make mistakes, and have their understanding challenged without feeling that every interaction is an examination.

I think this is a good place to stop for now. Thanks for the conversation, it has been great!

\speaker{AK} Thanks for chatting!

\section*{A logical map of the conversation}
\noindent\emph{The arrows below record dependencies suggested in the conversation, not formal logical implications. The Solved Axiom is a standing assumption throughout, while the Doom, Training, and Funding Premises are adverse predictions considered as part of the thought experiment. Several of the other axioms are deliberately provisional.}

\begin{center}
\begin{tikzpicture}[
  x=1cm,y=1cm,
  >={Latex[length=2.2mm]},
  every node/.style={font=\small,align=center},
  axiom/.style={draw=statementframe,rounded corners=2pt,fill=statementbg,minimum height=9mm,text width=29mm},
  theorem/.style={draw=akcolor,rounded corners=2pt,fill=white,minimum height=9mm,text width=32mm,line width=.7pt},
  final/.style={draw=tkcolor,rounded corners=2pt,fill=white,minimum height=10mm,text width=38mm,line width=.9pt},
  arr/.style={->,line width=.55pt,draw=stagegray}
]

\node[axiom] (solved) at (-5.7,3.2) {Solved Axiom};
\node[axiom] (epistemic) at (-2.1,3.2) {Epistemic Axiom};
\node[axiom] (learning) at (1.5,3.2) {Learning Axiom};
\node[axiom] (training) at (5,3.2) {Training Axiom};

\node[axiom] (value) at (-2.3,1.65) {Value Axiom};
\node[theorem] (expertise) at (1.5,1.65) {Expertise Demand Theorem};

\node[axiom] (dialectical) at (-4.7,0) {Dialectical Understanding Axiom};
\node[theorem] (trainingdemand) at (-1.0,0) {Training Demand Theorem};
\node[theorem] (hiring) at (3.2,0) {Hiring Value Theorem};

\node[final] (final) at (1.15,-1.65) {The Dialectical Theorem};

\draw[arr] (epistemic) -- (value);
\draw[arr] (value) -- (expertise);
\draw[arr] (learning) -- (expertise);
\draw[arr] (training) -- (expertise);
\draw[arr] (solved) -- (dialectical);
\draw[arr] (expertise) -- (trainingdemand);
\draw[arr] (expertise) -- (hiring);
\draw[arr] (dialectical) -- (final);
\draw[arr] (hiring) -- (final);

\end{tikzpicture}
\end{center}

\section*{Acknowledgements}

This conversation was conducted at the American Mathematical Society Office of Government Relations in Washington, D.C., and AK and TK both thank the AMS for its hospitality. AK and TK thank Chris Connell, Sara Maloni, Grant Sanderson, and Karen Saxe for constructive comments. Polishing, tightening, and adversarial reading of this conversation were provided by OpenAI's ChatGPT; the text was all human-written.

\printbibliography[title={References}]

\bigskip
\begin{center}
\small
\textbf{Thomas Koberda}\\
Department of Mathematics, University of Virginia\\
\href{mailto:thomas.koberda@virginia.edu}
     {thomas.koberda@virginia.edu}\\
\href{https://sites.google.com/view/koberdat}
     {Personal website}

\medskip

\textbf{Alex Kontorovich}\\
Department of Mathematics, Rutgers University\\
\href{mailto:alex.kontorovich@rutgers.edu}
     {alex.kontorovich@rutgers.edu}\\
\href{https://sites.math.rutgers.edu/~alexk/}
     {Personal website}
\end{center}
\end{document}